\documentclass[12pt,a4paper]{article}
\usepackage{hyperref}
\usepackage{geometry}
\usepackage{amsmath,amsthm,amscd,amsfonts,amssymb,enumerate}
\usepackage{tikz}
\usetikzlibrary{calc}
\usetikzlibrary{scopes}
\usetikzlibrary{backgrounds}
\usetikzlibrary{fit}
\usetikzlibrary{matrix}
\usetikzlibrary{positioning}

\newcommand{\drawBracket}[7]{
	\draw[line width=#7,rounded corners=#6]    (#1+#5*#4-#5*#2,#2-#5*#3+#5*#1) -- 
	(#1,#2)-- (0.5*#1+0.5*#3,0.5*#2+0.5*#4) -- (0.5*#1+0.5*#3-#5*#4+#5*#2,0.5*#2+0.5*#4+#5*#3-#5*#1);
	\draw[line width=#7,rounded corners=#6]
	(0.5*#1+0.5*#3-#5*#4+#5*#2,0.5*#2+0.5*#4+#5*#3-#5*#1)--
	(0.5*#1+0.5*#3,0.5*#2+0.5*#4) 
	-- (#3,#4)  --(#3-#5*#2+#5*#4,#4-#5*#3+#5*#1) ;
}%
\newtheorem{preproof}{{\bf Proof}\hspace{-.15cm}}

\newcommand{\qs}{\preceq _s}

\newcommand{\s}{\prec_s}
\newcommand{\ssq}{\subseteq}
\newcommand{\proofbegin}{\begin{preproof}\rm}
	\newcommand{\proofend}{\hfill{$\square$} \end{preproof}}
\newtheorem{thm}{Theorem}[section]
\newcommand{\theobegin}{\begin{thm} }
	\newcommand{\theoend}{\end{thm}}
\newtheorem{lem}[thm]{Lemma}
\newcommand{\lembegin}{\begin{lem} }
	\newcommand{\lemend}{\end{lem}}
\newtheorem{rem}[thm]{Remark }
\newcommand{\rembegin}{\begin{rem}\rm }
	\newcommand{\remend}{\end{rem}}
\newtheorem{prop}[thm]{Proposition}
\newcommand{\propbegin}{\begin{prop} }
	\newcommand{\propend}{\end{prop}}
\newtheorem{conj}[thm]{Conjecture}
\newcommand{\conjbegin}{\begin{conj}\rm }
	\newcommand{\conjend}{\end{conj}}
\newtheorem{cor}[thm]{Corollary}
\newcommand{\corbegin}{\begin{cor}\rm }
	\newcommand{\corend}{\end{cor}}
\newtheorem{question}[thm]{Question}
\newcommand{\questionbegin}{\begin{question}\rm }
	\newcommand{\questionend}{\end{question}}
\newtheorem{defin}[thm]{Definition}
\newcommand{\defbegin}{\begin{defin} }
	\newcommand{\defend}{\end{defin}}
\newtheorem{preexample}{{\bf Example}\hspace{-.15cm}}

\newcommand{\examplebegin}{\begin{preexample}\rm}
	\newcommand{\exampleend}{\end{preexample}}
\date{}
\title{
	\baselineskip = 0.8cm \vskip 1cm \bf The signless Laplacian Estrada index of tricyclic graphs}
\author{
	\bf R. Nasiri$^a$, H. R. Elahi$^a$, G. H. Fath-Tabar$^{b,}$\footnote{Corresponding author} , A. Gholami$^a$,\\ \bf T. Do\v sli\'c$^{ c}$\footnote{Email addresses: \href{mailto:R.Nasiri.82@gmail.com}{R.Nasiri.82@gmail.com} (R. Nasiri)
		\href{mailto:h.r.ellahi@gmail.com}{H.R.Ellahi@gmail.com} (H. R. Ellahi)
		\href{mailto:FathTabar@Kashanu.ac.ir}{FathTabar@Kashanu.ac.ir} (G. H. Fath-Tabar)
		\href{mailto:A.Gholami@Qom.ac.ir}{A.Gholami@Qom.ac.ir} (A. Gholami)
		\href{mailto:doslic@grad.hr}{doslic@grad.hr} (T. Do\v sli\'c)
	}\\
	\small  $^a$ Department of Mathematics, Faculty of  Sciences, University of Qom,\\
	\small Qom 37161-46611, I. R. Iran
	\\
	\small $^b$ Department of Pure Mathematics, Faculty of
	Mathematical Sciences,\\
	\small University of Kashan,  Kashan
	87317-51167, I. R. Iran\\
	\small $^c$ Faculty of Civil Engineering, University of Zagreb, Zagreb, Croatia }

\date{}
\begin{document}
\maketitle

\begin{abstract}
	\baselineskip=0.5cm
	\noindent 
	The signless Laplacian Estrada index of a graph $G$ is defined as $SLEE(G)=\sum^{n}_{i=1}e^{q_i}$ where $q_1, q_2, \ldots, q_n$ are the eigenvalues of the signless Laplacian matrix of $G$. 
	In this paper, we show that there are exactly two tricyclic graphs with the maximal signless Laplacian Estrada index.
	
	\vskip 3mm
	
	
	\vskip 3mm
	
\end{abstract}

\section{ Introduction}

Throughout this work, we are concerned with undirected simple graphs.
The vertex and edge sets of a graph $G$ are denoted by $V(G)$ and $E(G)$, 
respectively, and we assume $\vert V(G) \vert=n$ and $\vert E(G) \vert=m$. 
The \emph{cyclomatic number} of a graph $G$ is defined as $r(G) = m-n+c$, where
$c$ is the number of connected components of $G$. For connected graphs it means
$r(G) = m-n+1$. 
If $r(G) = 3$, i.e, $m=n+2$, then $G$ is called a \emph{tricyclic graph}.
The class of all tricyclic graphs on $n$ vertices is denoted by $\mathcal{J}_n$.

The adjacency matrix $\mathbf{A}=\mathbf{A}(G)=[a_{ij}]$ of $G$ is a matrix whose $(i,j)$-th entry is equal to 1 if vertices $i$ and $j$ are adjacent, and
0 otherwise. 
The set of all eigenvalues of $A(G)$ is the \emph{spectrum} of $G$, and the
largest eigenvalue of $G$ is called the \emph{spectral radius} of $G$. It is 
well-known that different graphs can have the same spectrum; the smallest
example is the so-called Saltire pair \cite{Dam}. Hence, the question arises
whether it is possible to uniquely reconstruct a graph from its spectrum. If
this is possible, i.e., if a graph $G$ is the only graph with a given spectrum,
we say that $G$ is \emph{determined by its spectrum} and abbreviate it as 
$G$ is DS. It is an open problem to determine the asymptotic fraction of DS
graphs; we refer the reader to \cite{Dam} for an in-deep treatment of this
problem.

Denote by $\mathbf{Q}=\mathbf{D}+\mathbf{A}$ the signless Laplacian matrix of $G$, where $\mathbf{D} = diag(d_1,d_2,\ldots,d_n)$ is the diagonal matrix of vertex degrees. The matrix $L = D - A$ is the usual Laplacian matrix of $G$.
The matrix $\mathbf{Q}$ is positive semi-definite, so the eigenvalues of $\mathbf{Q}$ can be ordered as $q^{}_{1}\geq  q^{}_{2}\geq \cdots\geq q^{}_{n}\geq 0$.
The largest eigenvalue of $\mathbf{Q}$ is called the \emph{signless Laplacian spectral radius} of graph and it is well known that this eigenvalue is simple and has a unique positive unit eigenvector.
In spectral graph theory, the problem of determining graphs with maximum spectral radius of $\mathbf{Q}$ is a prominent one on which many scholars have worked (e.g. \cite{Fan02,Tam}). 
Additional results about spectral properties of the signless Laplacian matrix have been reported in \cite{Abreu,Cvetkovic02,Cvetkovic03,Zhang}.
In \cite{Dam}, the authors present some evidence that the matrix $\mathbf{Q}$ 
might be better suited than the other graph matrices for studying various
graph properties.

The \emph{Estrada index} $EE(G)$ of a graph $G$ is defined as 
$$EE(G) = \sum _{i=1}^n e^{\lambda _j},$$
where $\lambda _j$ are the eigenvalues of $G$. It was introduced by Estrada in
2000 \cite{estrada} and has since found wide applications in chemistry, 
protein folding and study of complex networks. It was soon generalized also to 
the Laplacian case \cite{zhougutman}, and recently 
Ayyaswamy et al. \cite{Ayyaswamy01} defined it also for the signless Laplacian 
matrix. The \emph{signless Laplacian Estrada index} (hereafter denoted by 
$SLEE$) is defined as
\begin{equation}\label{eq:0}
	SLEE(G)=\sum^{n}_{i=1}e^{q_i}.
\end{equation}
Also, they specified bounds for $SLEE$ in terms of the number of vertices and edges.
Recently, Binthiya et al. \cite{Bin} established upper bound for $SLEE$ in 
terms of the vertex connectivity of a graph and defined the corresponding 
extremal graph.

Previously in \cite{Elahi01}, we characterized the unique graphs with maximum $SLEE$ among the set of all graphs with a given number of cut edges, pendent vertices, connectivity and edge connectivity. 
In this paper, we continue our study of $SLEE$ by determining the graphs 
with maximum $SLEE$ in the set of all tricyclic graphs. 
First and foremost, in section \ref{sec:2}, we state our main results in the form of two main theorems (\ref{theo:1} and \ref{theo:2}).
Afterwards, in section \ref{sec:3}, we provide some preliminary  definitions, 
notations and lemmas that we need to establish the validity of our main results.
Finally, in the last two sections we present our proofs of the theorems 
stated in  section \ref{sec:2}.

\section{Main results}\label{sec:2}
Our goal in this work is to prove the following theorems:
\theobegin\label{theo:1}
Let $G$ be a tricyclic graph on $n$ vertices having exactly $j$ simple cycles (i.e. a cycle in which no vertices are repeated), where $j\in \{3,4,6,7\}$.
If $G$ is an extremal graph with maximum $SLEE$, then $G\cong H_{j}^{n}$, where $H_{j}^{n}$ are as shown in Fig. 1, for  $j=3,4,6,7$.
\theoend
\begin{figure}[h]\label{fig:1}
	\center
	\begin{tikzpicture}
	\draw (-1.5,2) rectangle (10.5,-2.7);
	\begin{scope}[shift={(0, 0)}]
	\coordinate (0) at (0,0);
	\coordinate (1) at (1,0.5);
	\coordinate (2) at (1,-0.5);
	\coordinate (3) at (-1,0.5);
	\coordinate (4) at (-1,-0.5);
	\coordinate (5) at (-0.5,-1);
	\coordinate (6) at (0.5,-1);
	\draw  (0) -- (1) -- (2) --(0) -- (3) -- (4) -- (0) -- (5) -- (6) -- (0);
	\draw [ draw=white](-0.5,1.7) -- node[sloped]{$n-7$}(0.5,1.7) ;	
	\drawBracket{-0.7}{1.3}{0.8}{1.3}{0.05}{3pt}{0.5pt}
	\begin{scope}[shift={(0)}, x={(0:1cm)}, y={(90:1cm)}]
	\coordinate (x1) at (-0.5,1);	
	\coordinate (x2) at (-0.2,1);	
	\coordinate (x3) at (0.6,1);
	\draw [draw=white] (x2)-- node[sloped]{$\cdots$} (x3);
	\foreach \y in {x1, x2, x3}
	{
		\filldraw (\y) circle (2pt);
		\draw (\y) -- (0);	
	}			
	\end{scope}
	\begin{scope}[shift={(0,0.3)}]
	\draw [draw=white] (-1,-2)-- node[sloped]{$H^{n}_{3}$} (1,-2);
	\draw [draw=white] (-1,-2.5)-- node[sloped]{($j=3$, $n\geq 7$)} (1,-2.5);		
	\end{scope}	
	\foreach \y in {0,1,...,6}
	{
		\filldraw (\y) circle (2pt);	
	}	
	\end{scope}	
	\begin{scope}[shift={(3, 0)}]
	\coordinate (0) at (0,0);
	\coordinate (1) at (190:1.1);
	\coordinate (2) at (215:1.3);
	\coordinate (3) at (240:1.1);
	\coordinate (4) at (300:1.1);
	\coordinate (5) at (350:1.1);
	\draw  (0) -- (1) -- (2) -- (3) -- (0) -- (4) -- (5) -- (0)-- (2);
	\draw [ draw=white](-0.5,1.7) -- node[sloped]{$n-6$}(0.5,1.7) ;	
	\drawBracket{-0.7}{1.3}{0.8}{1.3}{0.05}{3pt}{0.5pt}
	\begin{scope}[shift={(0)}, x={(0:1cm)}, y={(90:1cm)}]		
	\coordinate (x1) at (-0.5,1);	
	\coordinate (x2) at (-0.2,1);	
	\coordinate (x3) at (0.6,1);
	\draw [draw=white] (x2)-- node[sloped]{$\cdots$} (x3);
	\foreach \y in {x1, x2, x3}
	{
		\draw (\y) -- (0);	
		\filldraw (\y) circle (2pt);
	}			
	\end{scope}
	\begin{scope}[shift={(0,0.3)}]
	\draw [draw=white] (-1,-2)-- node[sloped]{$H^{n}_{4}$} (1,-2);
	\draw [draw=white] (-1,-2.5)-- node[sloped]{($j=4$, $n\geq 6$)} (1,-2.5);			
	\end{scope}			\foreach \y in {0,1,...,5}
	{
		\filldraw (\y) circle (2pt);	
	}	
	\end{scope}	
	\begin{scope}[shift={(6, 0)}]
	\coordinate (0) at (0,0);
	\coordinate (1) at (-1,-0.5);
	\coordinate (2) at (0,-1);
	\coordinate (3) at (1,-0.5);
	\coordinate (4) at (0.5,-0.5);
	\draw  (0) -- (1) -- (2) -- (3)  -- (0) -- (2) -- (4) --(0);
	\draw [ draw=white](-0.5,1.7) -- node[sloped]{$n-5$}(0.5,1.7) ;	
	\drawBracket{-0.7}{1.3}{0.8}{1.3}{0.05}{3pt}{0.5pt}
	\begin{scope}[shift={(0)}, x={(0:1cm)}]		
	\coordinate (x1) at (-0.5,1);	
	\coordinate (x2) at (-0.2,1);	
	\coordinate (x3) at (0.6,1);
	\draw [draw=white] (x2)-- node[sloped]{$\cdots$} (x3);
	\foreach \y in {x1, x2, x3}
	{
		\draw (\y) -- (0);	
		\filldraw (\y) circle (2pt);
	}			
	\end{scope}
	\begin{scope}[shift={(0,0.3)}]
	\draw [draw=white] (-1,-2)-- node[sloped]{$H^{n}_{6}$} (1,-2);
	\draw [draw=white] (-1,-2.5)-- node[sloped]{($j=6$, $n\geq 5$)} (1,-2.5);
	\end{scope}	
	\foreach \y in {0,1,...,4}
	{
		\filldraw (\y) circle (2pt);	
	}	
	\end{scope}
	\begin{scope}[shift={(9, 0)}]
	\coordinate (0) at (0,0);
	\coordinate (1) at (-0.65,-1);
	\coordinate (2) at (0.65,-1);
	\coordinate (3) at (0,-0.65);
	\draw (0)--(1)--(2)--(3)--(1);
	\draw (2)--(0)--(3);
	\draw [ draw=white](-0.5,1.7) -- node[sloped]{$n-4$}(0.5,1.7) ;	
	\drawBracket{-0.7}{1.3}{0.8}{1.3}{0.05}{3pt}{0.5pt}
	\begin{scope}[shift={(0)}]		
	\coordinate (x1) at (-0.5,1);	
	\coordinate (x2) at (-0.2,1);	
	\coordinate (x3) at (0.6,1);
	\draw [draw=white] (x2)-- node[sloped]{$\cdots$} (x3);
	\foreach \y in {x1, x2, x3}
	{
		\draw (\y) -- (0);	
		\filldraw (\y) circle (2pt);
	}			
	\end{scope}
	\begin{scope}[shift={(0,0.3)}]
	\draw [draw=white] (-1,-2)-- node[sloped]{$H^{n}_{7}$} (1,-2);
	\draw [draw=white] (-1,-2.5)-- node[sloped]{($j=7$, $n\geq 4$)} (1,-2.5);
	\end{scope}
	\foreach \y in {0,1,...,3}		
	{
		\filldraw (\y) circle (2pt);	
	}	
	\end{scope}
	\end{tikzpicture}
	\vspace{3mm}
	\caption{The extremal tricyclic graphs with maximum $SLEE$ and the given number of simple cycles}
\end{figure}

\theobegin\label{theo:2}
If $G$ is a tricyclic graph on $n$ vertices, then $SLEE(G)\leq SLEE(H_{6}^{n})=SLEE(H_{7}^{n})$, with equality if and only if either $G\cong H_{6}^{n}$ or  $G\cong H_{7}^{n}$.
\theoend

In other words, theorem \ref{theo:1} shows the extremal unique graph with maximum $SLEE$ in the set of all $n$-vertex tricyclic graphs with the given number of simple cycles $j$, where $j\in\{3,4,6,7\}$.

Theorem \ref{theo:2} characterizes extremal tricyclic graphs with maximum $SLEE$ in the set of all tricyclic graphs on $n$ vertices, where $n\geq 5$.
Note that there are no tricyclic graphs on $n\leq 3$ vertices. 
Moreover, in the case where $n=4$, $H_{7}^{4}$ (which is the complete graph on $4$ vertices) is the unique (extremal) tricyclic graph.

There is an important point in the content of theorem \ref{theo:2}.
In fact, it consists of two propositions: 
First, of course, it proposes two graphs $H_{6}^{n}$ and $H_{7}^{n}$ to be the only possible extremal $n$-vertex tricyclic graphs with maximum $SLEE$;
second, it expresses that $SLEE(H_{6}^{n})=SLEE(H_{7}^{n})$. 
Actually, to prove the latter, we show (in section \ref{sec:5}) that the graphs  $H_{6}^{n}$ and $H_{7}^{n}$ have the equal signless Laplacian eigenvalues sequences.

\section{Preliminaries and lemmas}\label{sec:3}

In this section, we first state some definitions and notations used in our 
research and restate some results proved in references 
\cite{Cvetkovic01,Elahi01}.
Then, we prove an auxiliary lemma which is important for achieving the goals of this article.

Recall that the \emph{$k$-th signless Laplacian spectral moment} of a graph $G$ is denoted by $T^{}_k(G)$, and defined by $T^{}_k(G)=\sum^{n}_{i=1}q^{k}_{i}$.
Let $\mathbf{Q}_{}^{k}$ be the $k$-th power of signless Laplacian matrix of the graph $G$.
Then, by definition of $T^{}_{k}(G)$, one can easily see that $T^{}_{k}(G)$ is equal to the trace of matrix $\mathbf{Q}_{}^{k}$, i.e. $T^{}_{k}(G)=Tr(\mathbf{Q}^{k}_{})$.
Therefore, by Taylor expansion of exponential function $e^{q^{}_{i}}$ and the definition of $SLEE(G)$, we conclude 
\begin{equation}\label{eq:1}
	SLEE(G)=\sum^{}_{k\geq 0}\frac{T^{}_k(G)}{k!}.
\end{equation}
This equation leads us to the idea of using the notion of signless Laplacian spectral moments of graphs to compare their $SLEE$'s.
To exploit this idea, we need a notion which is very closely related to the signless Laplacian spectral moments of a graph.
The following definition and the proposition after that provide this suitable notion for us and state this close relation.
\defbegin\cite{Cvetkovic01}
A \emph{semi-edge walk} of length $k$ in graph $G$, is an alternating sequence 
$W=v^{}_1 e^{}_1 v^{}_2 e^{}_2 \cdots  v^{}_k e^{}_k v^{}_{k+1}$ of vertices $v^{}_1, v^{}_2, \dots , v^{}_k, v^{}_{k+1}$ and edges $e^{}_1, e^{}_2, \dots , e^{}_k$ such that  the vertices $v^{}_i$ and $v^{}_{i+1}$ are end-vertices (not necessarily distinct) of edge $e^{}_i$, for any $i=1, 2, \dots , k$.
If $v^{}_1=v^{}_{k+1}$, then we say $W$ is a \emph{closed semi-edge walk}.
\defend
\propbegin\cite{Cvetkovic01}\label{prop:1}
Let $G$ be a graph and $\mathbf{Q}$ be its signless Laplacian matrix.
The number of semi-edge walks of length $k$ starting at vertex $v$ and ending at vertex $u$ is equal to the $(v, u)$-th entry of the matrix $\mathbf{Q}^{k}_{}$.
\propend
As a consequence of the above proposition, it follows that the number of closed semi-edge walks of length $k$ in a graph $G$ is equal to $T^{}_k(G)$. 
Therefore, to calculate $T^{}_k(G)$ we may use the set of all (closed) semi-edge walks of length $k$ and its cardinality.

Let $G$ and $H$ be two graphs,  $x,y\in V(G)$, and $u,v\in V(H)$.
We denote by $SW^{}_k(G;x,y)$,  the set of all semi-edge walks, each of which is of length $k$ in $G$, starting at vertex $x$, and ending at vertex $y$. 
For convenience, we may denote  $SW^{}_k(G;x,x)$ by $SW^{}_k(G;x)$, and set $SW^{}_k(G)=\bigcup^{}_{x\in V(G)}SW^{}_k (G;x)$.
Given the above argument, we know that $T^{}_k(G)=|SW^{}_k(G)|$.
\\
We use the notation $(G;x,y)\qs(H;u,v)$ for, if  $ |SW^{}_k (G;x,y)|\leq|SW^{}_k(H;u,v)|$, for all $k\geq 0$.
Moreover, if  $(G;x,y)\qs(H;u,v)$,
and there exists some $k^{}_0$ such that $|SW^{}_{k^{}_0}(G;x,y)|<|SW^{}_{k^{}_0} (H;u,v)|$, then we write $(G;x,y)\s(H;u,v)$.

Now, we restate the following lemma, which is an excellent tool to compare the $SLEE$'s of two graphs, each of which  has a particular subgraph that are isomorphic.
\lembegin\cite{Elahi01}\label{lem:1}
Let $G$ be a graph and $v, u, w_1, w_2, \ldots , w_r\in V(G)$.
Suppose that
$E_v=\{e_1=vw_1, \ldots , e_r=vw_r\}$ and 
$E_u=\{e'_1=uw_1, \ldots ,  e'_r=uw_r\}$ where $e_i,e'_i\not\in E(G)$, for $i=1, 2, \ldots , r$.
Let $G_u=G+E_u$ and $G_v=G+E_v$.
If $(G;v)\s (G;u)$, and $(G;w_i,v)\qs (G;w_i,u)$ for each $i=1,2,\ldots,r$,
then $SLEE(G_v)<SLEE(G_u)$.
\lemend

To use the above lemma, we say that the graph $G^{}_{u}$ is obtained from $G^{}_{v}$ by transferring vertices  $w^{}_{1}, \dots, w^{}_{r}$ from  $N(v)$, the set of neighbors of $v$, to $N(u)$. 
In this situation, we call the vertices  $w^{}_{1}, \dots, w^{}_{r}$ the   \emph{transferred neighbors}, and the graph $G$ is named the  \emph{transfer route}.
Note that $G$ is a subgraph of both $G_u$ and $G_v$.

Although lemma \ref{lem:1} is an excellent tool, it has many conditions which have to be provided when we want to use it.
The following lemma enables us to discover a special case that provides such conditions.

\lembegin\label{lem:2}
Let $G$ be a graph and $u,v\in V(G)$. 
If $N(v)\ssq N(u)\cup \{u\}$, then $(G;v)\qs (G;u)$, and $(G;w,v)\qs (G;w,u)$ for each $w\in V(G)$. 
Moreover, if $deg_G(v)<deg_G(u)$, then $(G;v)\s (G;u)$, where $deg_G(v)$ is the degree of vertex $v$ in the graph $G$.
\lemend
\proofbegin
Let $k\geq 0$, and $W\in SW_k(G;v)$. 
We can decompose $W$ uniquely to $W_1W_2W_3$, such that $W_1$ and $W_2$ are as long as possible and consist of just the vertex $v$ and edges $vw$ where $w\in N(v)\setminus \{u\}$. 
Note that $W_2$ and $W_3$ are empty if $W$ does not contain any other vertex than $v$.
Let $W'_j$  be obtained from $W_j$, for $j=1,3$, by replacing the edge $v$ by $u$, and edges $vw$ by $uw$, where $w\in N(v)\setminus \{u\}$.
The map $f:SW_k(G;v)\to SW_k(G;u)$, defined by the rule $f(W_1W_2W_3)=W'_1W_2W'_3$, is injective.
Thus $(G;v)\qs (G;u)$.
\\
Similarly, by decomposing each semi-edge walk in $SW_k(G;w,v)$ and changing the end part of them, we conclude that $(G;w,v)\qs (G;w,u)$.
\\
Finally, since $deg_{G}^{}(v)=SW_1(G;v)$, the last part of the lemma is obvious.
\proofend
\section{The proof of theorem \ref{theo:1}}\label{sec:4}
In this section, we determine the unique $n$-vertex extremal tricyclic graph with maximum $SLEE$ which exactly has  $j$ simple cycles, for $j=3,4,6,7$.

The \emph{base} of a tricyclic graph $G$, denoted by $B(G)$, and defined to be the unique maximal subgraph of $G$, contains no pendent vertex.
Indeed, $B(G)$ is the unique minimal tricyclic subgraph of $G$, and $G$ can be obtained from $B(G)$ by planting some trees on it.

The following lemma is expressing the importance of the rule of the base of an extremal graph with maximum $SLEE$.

\lembegin\label{lem31}
Let $G$ be an extremal graph with maximum $SLEE$ over $\mathcal{J}_n$. Then 
each vertex in $G$ is either in $B(G)$ or a pendent vertex.
\lemend
\proofbegin
Let $T$ be a subgraph of $G$ with exactly one common vertex in $B(G)$, say $u$.
If $T$ is not a star with center vertex $u$, then there is a neighbor of $u$  in $T$, say  $v$, such that $deg_{G}^{}(v)>1$. 
Let   $G'$ be the graph obtained from $G$ by transferring all of the vertices in  $N(v)\setminus \{u\}$ from  $N(v)$ to $N(u)$, and $H$ be the transfer route graph.
Then, lemma \ref{lem:2} implies that $(H;v)\s(H;u)$, and therefore  $SLEE(G)<SLEE(G')$, by lemma \ref{lem:1}, which is a contradiction.
Hence, each subgraph of $G$ with just one common vertex $u$ in $B(G)$,  is a star with center vertex $u$, and this is the desired result.
\proofend
By \cite{Li01}, we know that there are 15 different bases of tricyclic graphs, and they can be classified  to four classes according to their number of simple cycles (as shown in Fig. 2):
\begin{align*}
	\mathcal{J}^3_n {}&= \{G\in\mathcal{J}^{}_n:B(G)\cong B^{3}_{i}: i\in\{1,2,\dots, 7\}\}\\
	\mathcal{J}^4_n {}&= \{G\in\mathcal{J}^{}_n:B(G)\cong B^{4}_{i}: i\in\{1,2,\dots, 4\}\}\\
	\mathcal{J}^6_n {}&= \{G\in\mathcal{J}^{}_n:B(G)\cong B^{6}_{i}: i\in\{1,2, 3\}\}\\
	\mathcal{J}^7_n {}&= \{G\in\mathcal{J}^{}_n:B(G)\cong B^{7}_{i}: i\in\{1\}\}
\end{align*}
With this classification, we surely have $\mathcal{J}^{}_n = \mathcal{J}^3_n\cup \mathcal{J}^4_n\cup \mathcal{J}^6_n\cup \mathcal{J}^7_n$.
Note that according to these notations, the aim of this section is to show that  $H_{j}^{n}$ is the unique extremal graph with maximum $SLEE$ among the members of $\mathcal{J}_{n}^{j}$, for any $j\in \{3,4,6,7\}$.
To reach this goal,
we need a suitable tool to compare the $SLEE$'s of graphs with the same type of bases.
Such a tool is provided in the following lemma.
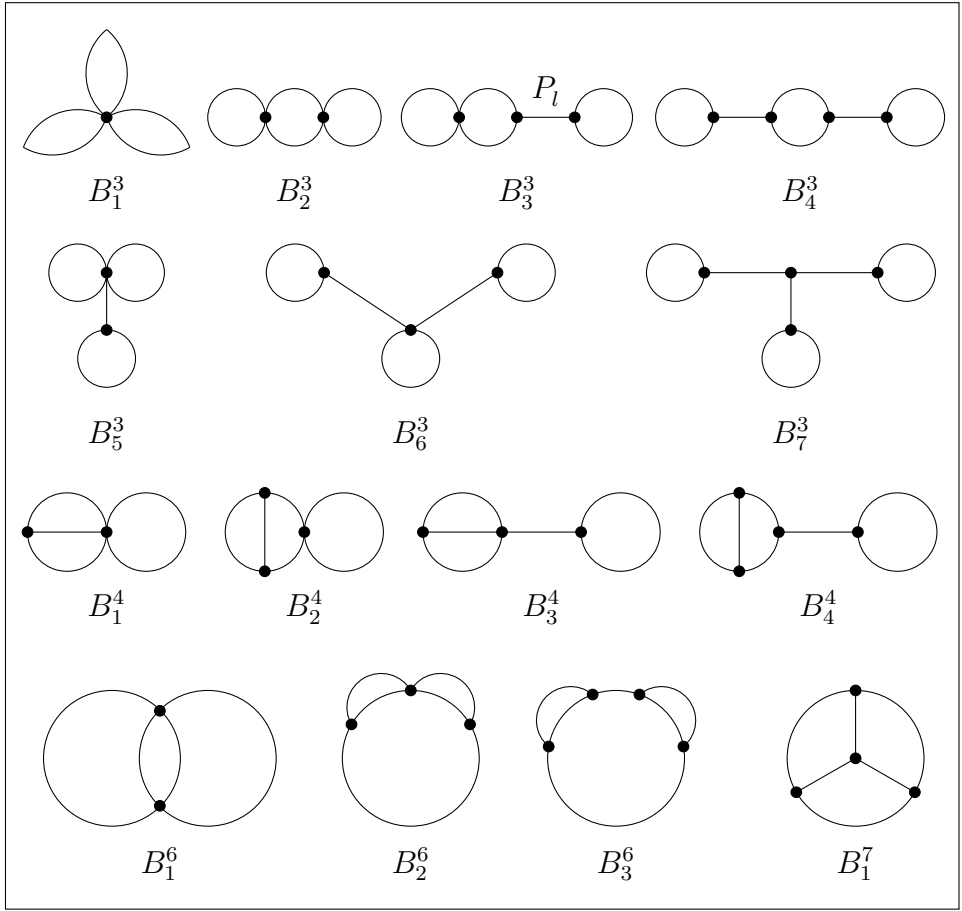
\begin{figure}[h]
	\center
	\begin{tikzpicture}
	\def\r{0.38}
	\draw (-3.5*\r,4*\r) rectangle (29.5*\r,-10.5);
	\begin{scope}
	
	\begin{scope}[shift={(0,0)}, scale=1]
	\draw [draw=white] (-1,-1)-- node[sloped]{$B^{3}_{1}$} (1,-1);
	\filldraw  (0,0) circle (2pt);
	\draw (0,0) arc (120:20:2*\r);
	\draw (0,0) arc (-160:-60:2*\r);
	\draw (0,0) arc (-50:50:2*\r);
	\draw (0,0) arc (230:130:2*\r);
	\draw (0,0) arc (60:160:2*\r);
	\draw (0,0) arc (-20:-120:2*\r);
	\end{scope}
	
	\begin{scope}[shift={(6.5*\r,0)}, scale=1]
	\draw [draw=white] (-1,-1)-- node[sloped]{$B^{3}_{2}$} (1,-1);
	\draw  (-2*\r,0) circle (\r);
	\draw  (0,0) circle (\r);
	\draw  (2*\r,0) circle (\r);
	\filldraw  (-\r,0) circle (2pt);
	\filldraw  (\r,0) circle (2pt);
	\end{scope}
	
	\begin{scope}[shift={(14.2*\r,0)}, scale=1]
	\draw [draw=white] (-1,-1)-- node[sloped]{$B^{3}_{3}$} (1,-1);
	\draw  (-3*\r,0) circle (\r);
	\draw  (-\r,0) circle (\r);
	\draw  (3*\r,0) circle (\r);
	\filldraw (-2*\r,0) circle (2pt);
	\filldraw (0,0) circle (2pt);
	\filldraw (2*\r,0) circle (2pt);
	\draw  (0,0)--node[anchor=south]{$P^{}_{l}$} (2*\r,0);
	\end{scope}
	
	\begin{scope}[shift={(24*\r,0)}, scale=1]
	\draw [draw=white] (-1,-1)-- node[sloped]{$B^{3}_{4}$} (1,-1);
	\draw (-4*\r,0) circle (\r);
	\draw  (0,0) circle (\r);
	\draw (4*\r,0) circle (\r);
	\draw  (-3*\r,0)--(-\r,0);
	\draw  (\r,0)--(3*\r,0);
	\filldraw  (-3*\r,0) circle (2pt);
	\filldraw  (-\r,0) circle (2pt);
	\filldraw  (\r,0) circle (2pt);
	\filldraw  (3*\r,0) circle (2pt);
	\end{scope}
	
	\end{scope}
	\begin{scope}
	
	\begin{scope}[shift={(0,-3.2)}, scale=1]
	\draw [draw=white] (-1,-1)-- node[sloped]{$B^{3}_{5}$} (1,-1);
	\draw (-\r,3*\r) circle (\r);
	\draw (0,0) circle (\r);
	\draw (\r,3*\r) circle (\r);
	\draw (0,3*\r) --(0,\r);
	
	\filldraw  (0,3*\r) circle (2pt);
	\filldraw  (0,\r) circle (2pt);
	\end{scope}
	
	\begin{scope}[shift={(4,-3.2)}, scale=1]
	\draw [draw=white] (-1,-1)-- node[sloped]{$B^{3}_{6}$} (1,-1);
	\coordinate (1) at (-3*\r,3*\r);
	\coordinate (2) at (0,\r);
	\coordinate (3) at (3*\r,3*\r);
	\draw (-4*\r,3*\r) circle (\r);
	\draw (0,0) circle (\r);
	\draw (4*\r,3*\r) circle (\r);
	\draw (1)-- (2)--(3);
	\filldraw (1) circle (2pt);
	\filldraw (2) circle (2pt);
	\filldraw (3) circle (2pt);
	\end{scope}
	
	\begin{scope}[shift={(9,-3.2)}, scale=1]
	\draw [draw=white] (-1,-1)-- node[sloped]{$B^{3}_{7}$} (1,-1);
	\coordinate (1) at (-3*\r,3*\r);
	\coordinate (2) at (0,\r);
	\coordinate (3) at (3*\r,3*\r);
	\coordinate (4) at (0,3*\r);
	\draw (-4*\r,3*\r) circle (\r);
	\draw (0,0) circle (\r);
	\draw (4*\r,3*\r) circle (\r);
	\draw (1) -- (3);
	\draw (4) -- (2);
	\filldraw (1) circle (2pt);
	\filldraw (2) circle (2pt);
	\filldraw (3) circle (2pt);
	\filldraw (4) circle (2pt);
	\end{scope}
	
	\end{scope}
	
	\def\r{0.52}
	
	\begin{scope}
	
	\begin{scope}[shift={(0,-5.5)}, scale=1]
	\draw [draw=white] (-1,-1) -- node[sloped]{$B^{4}_{1}$}(1,-1);
	\draw (-\r,0) circle (\r);
	\draw (\r,0) circle (\r);
	\draw (-2*\r,0)--(0,0);
	\filldraw (0,0) circle (2pt);
	\filldraw (-2*\r,0) circle (2pt);
	\end{scope}
	
	\begin{scope}[shift={(5*\r,-5.5)}, scale=1]
	\draw [draw=white] (-1,-1) -- node[sloped]{$B^{4}_{2}$}(1,-1);
	\draw (-\r,0) circle (\r);
	\draw (\r,0) circle (\r);
	\draw (-\r,\r)--(-\r,-\r);
	\filldraw  (0,0) circle (2pt);
	\filldraw  (-\r,\r) circle (2pt);
	\filldraw  (-\r,-\r) circle (2pt);
	\end{scope}
	
	\begin{scope}[shift={(11*\r,-5.5)}, scale=1]
	\draw [draw=white] (-1,-1) -- node[sloped]{$B^{4}_{3}$}(1,-1);
	\draw (-2*\r,0) circle (\r);
	\draw (2*\r,0) circle (\r);
	\draw (-3*\r,0)--(\r,0);
	\filldraw (\r,0) circle (2pt);
	\filldraw (-3*\r,0) circle (2pt);
	\filldraw (-\r,0) circle (2pt);
	\end{scope}
	
	\begin{scope}[shift={(18*\r,-5.5)}, scale=1]
	\draw [draw=white] (-1,-1) -- node[sloped]{$B^{4}_{4}$}(1,-1);
	\draw (-2*\r,0) circle (\r);
	\draw (2*\r,0) circle (\r);
	\draw (-2*\r,\r)--(-2*\r,-\r);
	\draw (-\r,0)--(\r,0);
	\filldraw  (\r,0) circle (2pt);
	\filldraw  (-\r,0) circle (2pt);
	\filldraw  (-2*\r,\r) circle (2pt);
	\filldraw  (-2*\r,-\r) circle (2pt);
	\end{scope}
	
	\end{scope}
	
	\def\r{0.9}
	
	\begin{scope}[shift={(0.4,0)}]
	
	\begin{scope}[shift={(0.3,-8.5)}, scale=1]
	\draw [draw=white] (-1,-\r-0.5)-- node[sloped]{$B^{6}_{1}$}(1,-\r-0.5);
	\draw  (-0.7*\r,0) circle (\r);
	\draw  (0.7*\r,0) circle (\r);
	\filldraw (0,0.7*\r) circle (2pt);
	\filldraw  (0,-0.7*\r) circle (2pt);
	\end{scope}
	
	\begin{scope}[shift={(4*\r,-8.5)}, scale=1]
	\draw  (0,0) circle (\r);
	\draw  (0,\r) arc (150:-30:0.5*\r);
	\draw  (0,\r) arc (30:210:0.5*\r);
	\filldraw  (0,\r) circle (2pt);
	\filldraw  (30:\r) circle (2pt);
	\filldraw  (150:\r) circle (2pt);
	\draw  [draw=white ] (-1,-\r-0.5)-- node[sloped]{$B^{6}_{2}$}(1,-\r-0.5);
	\end{scope}
	
	\begin{scope}[shift={(7*\r,-8.5)}, scale=1]
	\draw  [draw=white] (-1,-\r-0.5)-- node[sloped]{$B^{6}_{3}$}(1,-\r-0.5);
	\draw  (0,0) circle (\r);
	\draw (70:\r) arc (130:-50:0.5*\r);
	\draw  (110:\r) arc (50:230:0.5*\r);
	\filldraw (70:\r) circle (2pt);
	\filldraw  (110:\r) circle (2pt);
	\filldraw  (10:\r) circle (2pt);
	\filldraw  (170:\r) circle (2pt);
	\end{scope}
	
	\begin{scope}[shift={(10.5*\r,-8.5)}, scale=1]
	\draw  [draw=white] (-1,-\r-0.5)-- node[sloped]{$B^{7}_{1}$}(1,-\r-0.5);
	\draw  (0,0) circle (\r);
	\draw (90:\r) -- (0,0) -- (-30:\r);
	\draw (0,0) -- (210:\r);
	\filldraw (0,0) circle (2pt);
	\filldraw  (90:\r) circle (2pt);
	\filldraw  (-30:\r) circle (2pt);
	\filldraw  (210:\r) circle (2pt);
	\end{scope}
	
	\end{scope}
	\end{tikzpicture}
	\caption{A demonstration of all possible bases in $\mathcal{J}^{}_n$.}
\end{figure}
\lembegin
Let  $G$ be a tricyclic graph with $u,v\in V(G)$, where $e=uv\in E(G)$ and $N(u)\cap N(v)=\emptyset$.
If $G\in \mathcal{J}^{j}_n$, for some $j\in\{3,4,6,7\}$, then there exists a graph $G' \in \mathcal{J}^{j}_n$ such that $SLEE(G)<SLEE(G')$.
\lemend
\proofbegin
Let $G'$ be the graph obtained from $G$, by transferring all of vertices in $N(v)\setminus \{u\}$ from $N(v)$ to  $N(u)$, and $H$ be the transfer route graph.
By lemma \ref{lem:2} we have $(H;v)\s(H;u)$, which, according to the lemma \ref{lem:1}, implies that $SLEE(G)<SLEE(G')$.
On the other hand, since the aforesaid transferring does not change the number of neither simple cycles nor edges, we conclude that $G'\in \mathcal{J}^{j}_n$.
\proofend

\rembegin
By previous lemma, if the base of a tricyclic graph $G$ has a path which is not in a simple cycle (e.g. $B^{3}_{3}$ has the path $P^{}_{l}$, as shown in Fig. 2), then $G$ is not an extremal graph with maximum $SLEE$ through $\mathcal{J}^{j}_n$, for $j\in \{3,4\}$.
Moreover, if two successive vertices of 
a simple cycle of $G$, say the vertices $v^{}_{1}$ and $v^{}_{2}$ of the simple cycle $c^{}_{q}=v^{}_{1}v^{}_{2}\cdots v^{}_{q}v^{}_{1}$,  have no common neighbors (i.e. $N(v^{}_{1})\cap N(v^{}_{2})=\emptyset$), then $G$ is not an extremal graph with maximum $SLEE$ in $\mathcal{J}^{j}_n$, for $j\in \{3,4,6,7\}$.
Therefore, there are just 7 specific bases $A_{i}^{j}$ (as shown in Fig. 3), which are candidates to be the base of an extremal graph with maximum $SLEE$ among the members of $\mathcal{J}^{j}_n$, where $j\in \{3,4,6,7\}$ and $i\in\{1,2\}$.
\remend

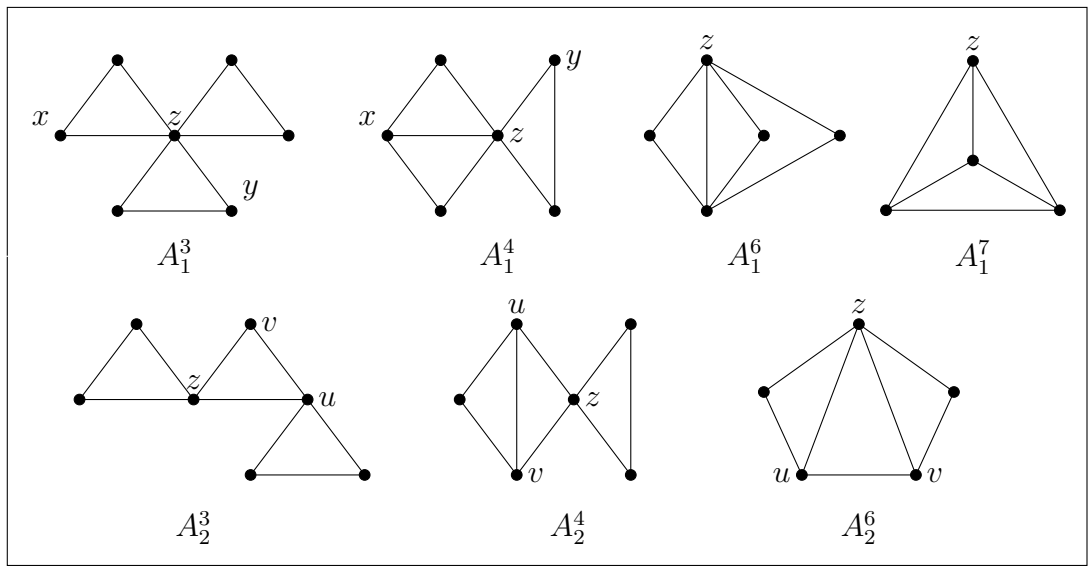
\begin{figure}[h]
	\center
	\begin{tikzpicture}
	\draw (-0.7,6.2) rectangle (13.5,-1.2);
	\begin{scope}[shift={(-9.25,0)}] 
	\filldraw (10.,5.5)circle (2pt);
	\filldraw (10.,3.5)circle (2pt);
	\filldraw (11.5,5.5)circle (2pt);
	\filldraw (11.5,3.5)circle (2pt)node[anchor=south west]{$y$};
	\filldraw (10.75,4.5)circle (2pt) node[anchor=south]{$z$};
	\filldraw (9.25,4.5)circle (2pt)node[anchor=south east]{$x$};
	\filldraw (12.25,4.5)circle (2pt);
	\draw (9.25,4.5) --(12.25,4.5)-- (11.5,5.5) --  (10.,3.5)  -- (11.5,3.5) --  (10.,5.5) --(9.25,4.5) ;
	\draw [draw=white] (8.5,2.9)-- node[sloped]{$A^{3}_{1}$} (13,2.9);
	\end{scope}
	\def\r{1.32}
	\begin{scope}[shift={(12,4.17)}] 
	\coordinate (0) at (0,0);
	\coordinate (1) at (90:\r);
	\coordinate (2) at (210:\r);
	\coordinate (3) at (-30:\r);
	\draw (1)--(0)--(2)--(1)--(3)--(2);
	\draw(0)--(3);
	\filldraw (0) circle (2pt);
	\filldraw (1) circle (2pt) node[anchor=south]{$z$};
	\filldraw (2) circle (2pt);
	\filldraw (3) circle (2pt);
	\draw [draw=white] (-1,-1.28)-- node[sloped]{$A^{7}_{1}$} (1,-1.28);
	\end{scope}
	
	\begin{scope}[shift={(6.5,0)}] 
	\filldraw (2.,5.5)circle (2pt) node[anchor=south]{$z$};
	\filldraw (2.,3.5)circle (2pt);
	\filldraw (1.25,4.5)circle (2pt);
	\filldraw (2.75,4.5)circle (2pt);
	\filldraw (3.75,4.5)circle (2pt);
	\draw  (2.,5.5) -- (3.75,4.5) -- (2.,3.5) -- (1.25,4.5) -- (2.,5.5);
	\draw (2.,5.5) -- (2.75,4.5)  -- (2.,3.5)-- (2.,5.5);
	\draw [draw=white] (.5,2.9)-- node[sloped]{$A^{6}_{1}$} (4.5,2.9);
	\end{scope}
	
	\begin{scope}[shift={(-1,0)}] 
	\filldraw (6.,5.5)circle (2pt);
	\filldraw (6.,3.5)circle (2pt);
	\filldraw (5.3,4.5)circle (2pt) node[anchor=south east]{$x$};
	\filldraw (6.75,4.5)circle (2pt) node[anchor=west]{$z$};
	\filldraw (7.5,5.5)circle (2pt)node[anchor=west ]{$y$};
	\filldraw (7.5,3.5)circle (2pt);
	\draw  (6.,5.5) -- (5.25,4.5) -- (6.,3.5)-- (6.75,4.5) -- (6.,5.5) ;
	\draw  (6.75,4.5) --(7.5,5.5) -- (7.5,3.5)--(6.75,4.5)--(5.25,4.5);
	\draw [draw=white] (5.,2.9)-- node[sloped]{$A^{4}_{1}$} (8.5,2.9);
	\end{scope}

	\begin{scope}[shift={(8,3.5)}] 
	\filldraw (1.25,1.6-4)circle (2pt);
	\filldraw (3.75,1.6-4)circle (2pt);
	\filldraw (1.75,0.5-4)circle (2pt)node[anchor= east ]{$u$};
	\filldraw (3.25,0.5-4)circle (2pt)node[anchor= west]{$v$};
	\filldraw (2.5,2.5-4)circle (2pt) node[anchor=south]{$z$};
	\draw  (2.5,2.5-4)-- (1.75,0.5-4);
	\draw  (2.5,2.5-4)-- (3.25,0.5-4) --  (1.75,0.5-4) -- (1.25,1.6-4) -- (2.5,2.5-4) -- (3.75,1.6-4) -- (3.25,0.5-4);
	\draw [draw=white] (1.25,-0.2-4)-- node[sloped]{$A^{6}_{2}$} (3.75,-0.2-4);
	\end{scope}
	
	\begin{scope}[shift={(0,3.5)}] 
	\filldraw (6.,2.5-4)circle (2pt)node[anchor=south ]{$u$};
	\filldraw (6.,0.5-4)circle (2pt)node[anchor=west]{$v$};
	\filldraw (5.25,1.5-4)circle (2pt);
	\filldraw (6.75,1.5-4)circle (2pt) node[anchor=west]{$z$};
	\filldraw (7.5,2.5-4)circle (2pt);
	\filldraw (7.5,0.5-4)circle (2pt);
	\draw  (6.,2.5-4) -- (5.25,1.5-4) -- (6.,0.5-4)-- (6.75,1.5-4) -- (6.,2.5-4) ;
	\draw  (6.75,1.5-4) --(7.5,2.5-4) -- (7.5,0.5-4)--(6.75,1.5-4);
	\draw (6.,2.5-4)--(6.,0.5-4);
	\draw [draw=white] (5.25,-0.2-4)-- node[sloped]{$A^{4}_{2}$} (8.1,-0.2-4);
	\end{scope}
	
	\begin{scope}[shift={(-9,3.5)}] 
	\filldraw (10,2.5-4)circle (2pt);
	\filldraw (11.5,0.5-4)circle (2pt);
	\filldraw (11.5,2.5-4)circle (2pt)node[anchor=west]{$v$};
	\filldraw (13,0.5-4)circle (2pt);
	\filldraw (10.75,1.5-4)circle (2pt) node[anchor=south]{$z$};
	\filldraw (9.25,1.5-4)circle (2pt);
	\filldraw (12.25,1.5-4)circle (2pt)node[anchor=west]{$u$};
	\draw(12.25,1.5-4) --(9.25,1.5-4) --(10,2.5-4) --  (10.75,1.5-4);
	\draw (12.25,1.5-4)-- (11.5,0.5-4)-- (13,0.5-4)--  (11.5,2.5-4)--(10.75,1.5-4) ;
	\draw [draw=white] (8.5,-0.2-4)-- node[sloped]{$A^{3}_{2}$} (13,-0.2-4);
	\end{scope}
	\end{tikzpicture}
	\caption{Bases $A_{i}^{j}$ and $A_{1}^{7}$, for $i\in \{1,2\}$ and $j\in \{3,4,6\}$.}
\end{figure}

\lembegin\label{lem00014}
If $G$ be an extremal graph with maximum $SLEE$ in $\mathcal{J}^j_n$, for $j=3,4,6,7$, then $B(G)\cong A^j_1$.
\lemend
\proofbegin
The case $j=7$ is obvious.
Let $j=3$, and $G$ be an extremal graph with maximum $SLEE$ in $\mathcal{J}^{3}_n$, such that $B(G)\cong A^{3}_{2}$.
Suppose that $G'$ be the graph obtaining from $G$ by transferring all of vertices in $N(u)\setminus\{z,v\}$  from $N(u)$ to $N(z)$, and $H$ be the transfer route graph.
Note that $B(G')\cong A^{3}_{1}$.
Since  $N^{}_{H}(u)\ssq N^{}_{H}(z)\cup \{z\}$, lemma \ref{lem:2} implies that $(H;u)\s(H;z)$.
Thus, by lemma \ref{lem:1}, $SLEE(G)<SLEE(G')$, which is a contradiction.
Therefore, if $G$ is an extremal graph with maximum $SLEE$ in $\mathcal{J}^{3}_n$, then $B(G)\cong A^{3}_{1}$.

For the case $i=4$ (respectively, $i=6$), the result follows by a similar method used above and transferring all of vertices in $N(u)\setminus\{v,z\}$ from $N(u)$ to $N(z)$ (respectively, $N(v)$), where vertices $v,u$ and $z$ are shown in Fig. 3.
\proofend

Now, we are ready to prove our first main theorem, which states that $H_{j}^{n}$ is, up to isomorphism, the unique extremal tricyclic graph, such that it has maximum $SLEE$ with the given number of simple cycles $j$, for $j\in\{3,4,6,7\}$.
\\

{\bf Proof of theorem \ref{theo:1}}
Let $j\in \{3,4,6,7\}$ and  $G$ be an extremal graph with maximum $SLEE$ through $\mathcal{J}^j_n$.
By previous lemmas, $G$ is obtained by attaching some pendent vertices to some vertices of $A^{j}_{1}$.
Let $x$ be a vertex of $A^{j}_{1}$, where $x\neq z$ and it has some pendent neighbors (the vertex $z$ is shown in Fig. 3).
Further, let $N^{np}_{}(x)$ be the set of all non-pendent neighbors of $x$.
Since $N^{np}_{}(x)\ssq N(z)\cup\{z\}$, by transferring pendent neighbors of 
$x$ from $N(x)$ to $N(z)$ and applying lemma \ref{lem:2}, we get a graph $G'$ such that $SLEE(G)<SLEE(G')$, which is a contradiction. 
Therefore, all of  $n-|V(A^j_1)|$  pendent vertices of $G$ are attached to $z$.
It means that $G$ is isomorphic to $H_{j}^{n}$
\hfill{$\square$}

\section{The proof of theorem \ref{theo:2}}\label{sec:5}

In this section, we prove the first proposition of theorem \ref{theo:2} in the next lemma. 
afterward, we represent the proof of the second proposition of the theorem to close the section.

\lembegin
Let $G$ be an extremal graph with the maximum $SLEE$  through the graphs in 
$\mathcal{J}_n$. Then $B(G)\cong A^j_1$, where $j\in\{6,7\}$.
\lemend
\proofbegin
Let $G$ be an extremal graph with maximum $SLEE$ in $\mathcal{J}_n$.
By lemma \ref{lem00014}, it is enough to show that $B(G)$ is isomorphic to neither $A^{3}_{1}$ nor $A^{4}_{1}$.

Let $B(G)\cong A^{3}_{1}$ or $A^{4}_{1}$ (as shown in Fig. 3).
Suppose that $G'$ be the graph obtained from $G$ by transferring all of 
vertices in $N(y)\setminus \{z\}$ from $N(y)$ to $N(x)$, and let 
$H$ be the transfer route graph.
Since $N^{}_{H}(y)\subset N^{}_{H}(x)$, we have $(H;y)\s (H;x)$ and  $(H;w,y)\qs (H;w,x)$, by lemma \ref{lem:2}, for each $w\in N(y)\setminus \{z\}$.
Thus by lemma \ref{lem:1}, we get $SLEE(G)<SLEE(G')$, which is a contradiction.
Therefore $B(G)\cong A^{6}_{1}$ or $A^{7}_{1}$.
\proofend

{\bf Proof of theorem \ref{theo:2}}
For the first claim of the theorem, we turn your attention to the fact that theorem \ref{theo:1} and the previous lemma guarantee that only graphs $H^{n}_{6}$ and $H^{n}_{7}$ are candidates for being extremal graphs with maximum $SLEE$ through $\mathcal{J}^{}_n$.

In the following, we show that these two graphs have simultaneous maximum $SLEE$ value through the members of  $\mathcal{J}^{}_n$.
To do this, it is enough to show that the characteristic polynomials of  $\mathbf{Q}(H^{n}_{6})$ and $\mathbf{Q}(H^{n}_{7})$ are equal.
In fact, this means that the eigenvalues sequences of these matrices  are equal.
Therefore, by equation \ref{eq:0},  $SLEE(H^{n}_{6})=SLEE(H^{n}_{7})$ as desired.

Let us use labels for vertices of  $H^{n}_{6}$ and $H^{n}_{7}$ as shown in Fig. 4.
\begin{figure}[h]
	\center
	\begin{tikzpicture}
	\draw (3.5,-1.2) rectangle (15,4);
	\foreach \xshift/\yshift in {11/0}
	{
		\filldraw (\xshift+0,\yshift+0) circle (2pt) node [anchor=south east]{$v_2$};
		\filldraw (\xshift+2.5,\yshift+0) circle (2pt) node[anchor=south west]{$v_3$};
		\filldraw (\xshift+1.25,\yshift+0.75) circle (2pt) node[anchor=west]{$v_4$};
		\filldraw (\xshift+1.25,\yshift+2.25) circle (2pt) node[anchor=east]{$v_1$};
		\filldraw (\xshift+0,\yshift+3.25) circle (2pt) node [anchor=south]{$v_5$};
		
		\draw (\xshift+0,\yshift+3.25)--(\xshift+1.25,\yshift+2.25) --(\xshift+0,\yshift+0) -- (\xshift+2.5,\yshift+0)  -- (\xshift+1.25,\yshift+2.25)-- (\xshift+1.25,\yshift+0.75)-- (\xshift+2.5,\yshift+0);
		\draw (\xshift+0,\yshift+0)-- (\xshift+1.25,\yshift+0.75) ;
		\draw [draw=white] (\xshift+0,\yshift-0.75)-- node[sloped]{$H^{n}_{7}$} (\xshift+2.5,\yshift-0.75);
		
	}
	
	\foreach \x/\y/\l in {12.25/2.25/v,6/2.25/u}
	{
		\filldraw (\x-0.625,\y+1) circle (2pt) node [anchor=south]{$\l_6$};
		\filldraw (\x+1.25,\y+1) circle (2pt) node [anchor=south]{$\l_n$};
		\filldraw (\x,\y+1) circle (2pt) node [anchor=south]{$\l_7$};
		\draw  (\x-0.625,\y+1)--(\x,\y)-- (\x+1.25,\y+1);
		\draw  (\x,\y+1)--(\x,\y);
		\foreach \z in {0.8,0.4,0.6}
		\filldraw (\x+\z,\y+1) circle (0.7pt);
	}
	
	\filldraw (6,2.25) circle (2pt) node [anchor=east]{$u_1$};
	\filldraw (6,0) circle (2pt) node [anchor=east]{$u_2$};
	\filldraw (5,1.125) circle (2pt) node [anchor=east]{$u_3$};
	\filldraw (7,1.125) circle (2pt) node [anchor=east]{$u_4$};
	\filldraw (8,1.125) circle (2pt) node [anchor=west]{$u_5$};
	\draw (8,1.125)--(6,2.25) --(5,1.125) --(6,0) -- (7,1.125)--(6,2.25) --(6,0) --(8,1.125) ;
	\draw [draw=white] (4,-0.75)-- node[sloped]{$H^{n}_{6}$} (8,-0.75);
	\end{tikzpicture}
	\caption{The two $n$-vertex tricyclic graphs which are simultaneous extremal with maximum $SLEE$ value.}
\end{figure}
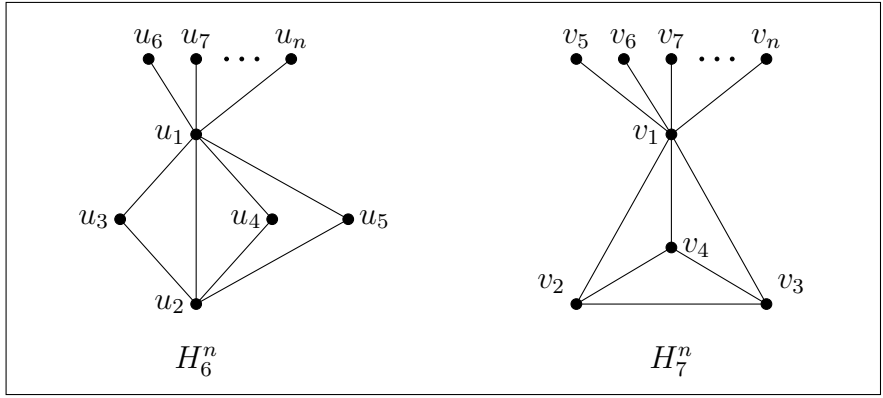

Let $j\in\{6,7\}$ and $\mathbf{S}^{}_{j}$  be the sub-matrix of  $\mathbf{Q}(H_{j}^{5})$, obtained by removing the first row and column:

$$\mathbf{S}_{6}^{}=\begin{aligned}
\begin{tikzpicture}
\matrix [matrix of math nodes,right delimiter=\rbrack,left delimiter= 
\lbrack] (m)
{
	4 & 1 & 1 & 1 \\
	1 & 2 & 0 & 0 \\
	1 & 0 & 2 & 0 \\
	1 & 0 & 0 & 2 \\
};
\end{tikzpicture}
\end{aligned}
\qquad
\mathbf{S}_{7}^{}=\begin{aligned}
\begin{tikzpicture}
\matrix [matrix of math nodes,right delimiter=\rbrack,left delimiter= 
\lbrack] (m)
{
	3 & 1 & 1 & 0\\
	1 & 3 & 1 & 0\\
	1 & 1 & 3 & 0\\
	0 & 0 & 0 & 1\\
};
\end{tikzpicture}
\end{aligned}
$$
With these notations, we have
$$\mathbf{Q}(H_{j}^{n}) - x\,\mathbf{I}_{n}^{}=\begin{aligned}
\begin{tikzpicture}
\matrix [matrix of math nodes,right delimiter=\rbrack,left delimiter= 
\lbrack] (m)
{
	4-x\, & 1\, & 1\, & 1\, & 1\, & 1 & \dots & 1\\
	1 &  . & .  &  . &  . & 0 & \dots & 0\\
	1 &  . & .  &  . &  . & 0 & \dots & 0\\
	1 &  . & .  &  . &  . & 0 & \dots & 0\\
	1 &  . & .  &  . &  . & 0 & \dots & 0\\
	1 & 0 & 0 & 0 & 0 & \, 1-x &  \dots     & 0\\
	\vdots & \vdots & \vdots & \vdots & \vdots &  \vdots   & \ddots & \vdots\\
	1 & 0 & 0 & 0 & 0 & 0 & \dots & \, 1-x\\
};
\def\d{0.03}
\coordinate (1) at (m-2-2.north west);
\coordinate (2) at (m-5-5.south east);
\path [fill=gray!30] (2) rectangle (1)+(\d,\d);
\draw (2) rectangle node{$\mathbf{S}_{j}^{} - x\,\mathbf{I}_{n}^{}$} (1)+(\d,\d);
\end{tikzpicture}
\end{aligned}
$$
Where $\mathbf{I}_{n}$ is the $n\times n$ identity matrix.
Note that the characteristic polynomial of any $n\times n$-matrix $\mathbf{M}$ is defined to be $det(\mathbf{M} - x\,\mathbf{I}_{n}^{})$.

Now, let $n\geq 6$.
We evaluate the $det(\mathbf{Q}(H_{j}^{n}) - x\,\mathbf{I}_{n}^{})$ by computing
the cofactor expansion along the last row and then computing the subsequent
cofactor expansions along the last column, and then we get the following formula for both $j=6$ and $7$:
$$det(\mathbf{Q}(H_{j}^{n}) - x\,\mathbf{I}_{n})=
(1-x)_{}^{n-6}\,\, det(\mathbf{S}_{j}^{} - x\,\mathbf{I}_n) + (1-x)\, det(\mathbf{Q}(H_{j}^{n-1}) - x\,\mathbf{I}_{n-1})
$$
Also, for $n=5$, by a few algebraic calculations, one can easily show that
$$det(\mathbf{Q}(H_{6}^{5}) - x\,\mathbf{I}_{5}^{})=-x^5+14\,x^4-69\,x^3+152\,x^2-148\,x+48=det(\mathbf{Q}(H_{7}^{5}) - x\,\mathbf{I}_{5}^{})$$
and moreover
$det(\mathbf{S}_{6} - x\,\mathbf{I}_{5}^{})=
x^4-10\,x^3+33\,x^2-44\,x+20
=det(\mathbf{S}_{7} - x\,\mathbf{I}_{5}^{})$.
Therefore, the result holds by induction on $n$ where $n\geq 5$ and using the last three formulas.
\hfill{$\square$}

\section{Concluding remarks}

In this paper we have studied the signless Laplacian Estrada index and 
determined its maximum values and corresponding graphs over the class of all
tricyclic graphs on a given number of vertices. It would be of interest to
study its behavior also on other classes of connected graphs with simple
connectivity patterns and cycle structure. Among the most promising 
candidates are linear polymers and cactus graphs, as well as the graphs that
arise from simpler components via binary operations of splice and link. 
\section*{Acknowledgments.} 

The third author has been financially supported by University of Kashan (Grant No. 504631/7).



\end{document}